# A Note On Testing Of Hypothesis

Rajesh Singh
School of Statistics, D. A.V.V., Indore (M.P.), India
rsinghstat@yahoo.co.in

Jayant Singh
Department of Statistics
Rajasthan University, Jaipur, India
Jayantsingh47@rediffmail.com

Florentin Smarandache
Chair of Department of Mathematics, University of New Mexico, Gallup, USA
fsmarandache@yahoo.com

Abstract : In this paper problem of testing of hypothesis is discussed when the samples have been drawn from normal distribution. The study of hypothesis testing is also extended to Baye's set up.

## 1. Introduction

Let the random variable (r.v.) X have a normal distribution $N(\theta,\ \sigma^2)$, $\sigma^2$ is assumed to be known. The hypothesis $H_0 : \theta = \theta_0$ against $H_1 : \theta = \theta_1$, $\theta_1 > \theta_0$ is to be tested. Let $X_1$, $X_2$, …, $X_n$ be a random sample from $N(\theta,\ \sigma^2)$ population. Let $\overline{X}\,(=\frac{1}{n}\sum_{i=1}^{n} X_i)$ be the sample mean.

By Neyman – Pearson lemma the most powerful test rejects $H_0$ at $\alpha$% level of significance,

if $\dfrac{\sqrt{n}\left(\overline{X}-\theta_o\right)}{\sigma} \geq d_\alpha$, where $d_\alpha$ is such that

$$\int_{d_\alpha}^{\infty} \frac{1}{\sqrt{2\pi}} e^{-\frac{Z^2}{2}} dZ = \alpha$$

If the sample is such that $H_0$ is rejected then will it imply that $H_1$ will be accepted?

In general this will not be true for all values of $\theta_1$, but will be true for some specific value of $\theta_1$ i.e., when $\theta_1$ is at a specific distance from $\theta_0$.

$H_0$ is rejected if $\dfrac{\sqrt{n}\left(\overline{X}-\theta_o\right)}{\sigma} \geq d_\alpha$

$$\text{i.e. } \overline{X} \geq \theta_0 + d_\alpha \frac{\sigma}{\sqrt{n}} \qquad (1)$$

Similarly the Most Powerful Test will accept $H_1$ against $H_0$

if $\dfrac{\sqrt{n}\left(\overline{X}-\theta_1\right)}{\sigma} \geq -d_\alpha$

i.e. $\quad \overline{X} \geq \theta_1 - d_\alpha \frac{\sigma}{\sqrt{n}}$ (2)

Rejecting $H_0$ will mean accepting $H_1$

if $\quad (1) \Rightarrow (2)$

i.e. $\quad \overline{X} \geq \theta_0 + d_\alpha \frac{\sigma}{\sqrt{n}} \quad \Rightarrow \quad \overline{X} \geq \theta_1 - d_\alpha \frac{\sigma}{\sqrt{n}}$

i.e. $\quad \theta_1 - d_\alpha \frac{\sigma}{\sqrt{n}} \leq \theta_0 + d_\alpha \frac{\sigma}{\sqrt{n}}$ (3)

Similarly accepting $H_1$ will mean rejecting $H_0$

if $\quad (2) \Rightarrow (1)$

i.e. $\quad \theta_0 + d_\alpha \frac{\sigma}{\sqrt{n}} \leq \theta_1 - d_\alpha \frac{\sigma}{\sqrt{n}}$ (4)

From (3) and (4) we have

$$\theta_0 + d_\alpha \frac{\sigma}{\sqrt{n}} = \theta_1 - d_\alpha \frac{\sigma}{\sqrt{n}}$$

i.e. $\quad \theta_1 - \theta_0 = 2\, d_\alpha \frac{\sigma}{\sqrt{n}}$ (5)

Thus $d_\alpha \frac{\sigma}{\sqrt{n}} = \frac{\theta_1 - \theta_0}{2}$ and $\theta_1 = \theta_0 + 2 d_\alpha \frac{\sigma}{\sqrt{n}}$.

From (1) Reject $H_0$ if $\overline{X} > \theta_0 + \frac{\theta_1 - \theta_0}{2} = \frac{\theta_0 + \theta_1}{2}$

and from (2) Accept $H_1$ if $\overline{X} > \theta_1 - \frac{\theta_1 - \theta_0}{2} = \frac{\theta_0 + \theta_1}{2}$

Thus rejecting $H_0$ will mean accepting $H_1$

when $\quad \overline{X} > \frac{\theta_0 + \theta_1}{2}$.

From (5) this will be true only when $\theta_1 = \theta_0 + 2\, d_\alpha \frac{\sigma}{\sqrt{n}}$. For other values of $\theta_1 \neq \theta_0 + 2\, d_\alpha \dfrac{\sigma}{\sqrt{n}}$ rejecting $H_0$ will not mean accepting $H_1$.

It is therefore, recommended that instead of testing $H_0 : \theta = \theta_0$ against $H_1 : \theta = \theta_1,\ \theta_1 > \theta_0$, it is more appropriate to test $H_0 : \theta = \theta_0$ against $H_1 : \theta > \theta_0$. In this situation rejecting $H_0$ will mean $\theta > \theta_0$ and is not equal to some given value $\theta_1$.

But in Baye's setup rejecting $H_0$ means accepting $H_1$ whatever may be $\theta_0$ and $\theta_1$. In this set up the level of significance is not a preassigned constant, but depends on $\theta_0$, $\theta_1$, $\sigma^2$ and n.

Consider (0,1) loss function and equal prior probabilities ½ for $\theta_0$ and $\theta_1$. The Baye's test rejects $H_0$ (accept $H_1$)

$$\text{if} \quad \bar{X} > \frac{\theta_0 + \theta_1}{2}$$

and accepts $H_0$ (rejects $H_1$)

$$\text{if} \quad \bar{X} < \frac{\theta_0 + \theta_1}{2}.$$

[See Rohatagi p.463, Example 2]

The level of significance is given by

$$P_{H_0}\ [\bar{X} > \frac{\theta_0 + \theta_1}{2}] = P_{H_0}[\ \frac{(\bar{X} - \theta_0)\sqrt{n}}{\sigma} > \frac{(\theta_1 - \theta_0)\sqrt{n}}{2\sigma}\ ]$$

$$= 1 - \Phi\left(\frac{\sqrt{n}(\theta_1 - \theta_0}{2\sigma}\right)$$

where $\Phi(t) = \int_{-\infty}^{t} \frac{1}{\sqrt{2\pi}} e^{-\frac{Z^2}{2}} dZ$.

Thus the level of significance depends on $\theta_0$, $\theta_1$, $\sigma^2$ and n.

Acknowledgement : Author’s are thankful to Prof. Jokhan Singh for suggesting this problem.